\documentclass[12pt]{article}
\usepackage[reqno]{amsmath}
\usepackage{amssymb, theorem, enumerate}
\usepackage{amsfonts}

\newcommand\al{\alpha}
\newcommand\be{\beta}
\newcommand\de{\delta}

\newcommand\cF{{\mathcal F}}

\newcommand\fld{{\mathbb F}}

\newcommand\diff{\mathbin{\mkern-1.5mu\setminus\mkern-1.5mu}}
\newcommand\gbin[2]{\genfrac{[}{]}{0pt}{}{#1}{#2}}
\newcommand\inv{^{-1}}

\newcommand\seq[3]{#1_{#2},\ldots,#1_{#3}}

\expandafter\ifx\csname pdfoptionalwaysusepdfpagebox\endcsname\relax\else
\fi

\theoremstyle{change}
{\theorembodyfont{\slshape}
\newtheorem{theorem}{Theorem.}[section]
\newtheorem{lemma}[theorem]{Lemma.}
}
\theorembodyfont{\rmfamily}

\newcommand\lref[1]{Lemma~\ref{lem:#1}}

\newcommand\cref[1]{Corollary~\ref{cor:#1}}
\newcommand\sref[1]{Section~\ref{sec:#1}}

\def\proof{\noindent{{\sl Proof. }}}
\def\sqr#1#2{{\vbox{\hrule height.#2pt
    \hbox{\vrule width.#2pt height#1pt \kern#1pt
        \vrule width.#2pt}\hrule height.#2pt}}}
\def\eqed{\sqr53}
\def\qed{%
    \ifmmode\eqno\eqed
    \else\nobreak\ \hfill\eqed\medbreak\fi}

\newcommand\kne[2]{K_{#1:#2}}
\newcommand\knvk{\kne{v}{k}}
\newcommand\qkne[2]{qK_{#1:#2}}
\newcommand\qkvk{\qkne{v}{k}}

\title{Colouring Lines in Projective Space}
\author{Ameera Chowdhury\\
    Mathematics\\
    Caltech\\
    Pasadena, CA 91125, U.S.A.
\and
Chris Godsil\\
    Combinatorics \& Optimisation\\
    University of Waterloo\\
    Waterloo, Ont. N2L 3G1, Canada
\and
Gordon Royle\\
    Computer Science \& Software Engineering\\
    University of Western Australia\\
    Crawley, WA 6009, Australia}

\begin{document}
\maketitle

\begin{abstract}
Let $V$ be a vector space of dimension $v$ over a field of order $q$.  The $q$-Kneser graph
has the $k$-dimensional subspaces of $V$ as its vertices, where two subspaces $\al$ and $\be$
are adjacent if and only if $\al\cap\be$ is the zero subspace.  This paper is motivated by the problem of determining the chromatic numbers of these graphs.  This problem is trivial when $k=1$ (and the graphs are complete) or when $v<2k$ (and the graphs are empty).  We establish some basic theory in the general case.   Then specializing to the case $k=2$, we show that the chromatic number is
$q^2+q$ when $v=4$ and $(q^{v-1}-1)/(q-1)$ when $v > 4$.  In both cases we characterise the
minimal colourings.
\end{abstract}

\section{Introduction}

The \textsl{Kneser graph} $K_{v:k}$ has the subsets of size $k$ from a fixed set of size $v$ as its vertices, with two $k$-subsets adjacent if and only if they are disjoint as sets.  The Kneser
graphs play an important role in work on graph homomorphisms and graph colouring.
In this paper we are concerned with a generalisation of these graphs, which we call
\textsl{$q$-Kneser graphs}.  We construct these as follows.  Let $\fld$ be a finite field of order $q$.  The vertices of the $q$-Kneser graph $\qkne{v}k$ are the $k$-dimensional subspaces of a vector space of dimension $v$ over $\fld$; two $k$-subspaces are adjacent if their intersection
is the zero subspace.

Our work in this paper is concerned with determining the chromatic numbers of the $q$-Kneser graphs.  The graphs $\qkne{v}1$ are complete graphs and if $v<2k$ then $\qkne{v}k$ is an empty graph, and there is nothing we need say about these cases.  We summarise our main results.

We show that if $v>2k$, then
\[
\chi(\qkne{v}k) \le\frac{q^{v-k+1}-1}{q-1}
\]
and if $v=2k$, then
\[
\chi(\qkne{2k}k) \le q^k+q^{k-1}.
\]
Naturally these bounds are derived by giving explicit colourings.  We prove that the stated bounds
are tight when $k=2$, where we can also characterise the minimal colourings.  When $v\ge5$ these
are essentially unique, but when $v=4$ there are a number of colourings.

\medbreak
We now explain why these colouring questions are interesting.  We first recall what is known about the ordinary Kneser graphs $\knvk$.  It is easy to find a colouring of $\knvk$ with $v-2k+2$ colours as follows: If $\al$ is a $k$-subset and the largest element of $\al$ is greater than $2k$, define this element to be the colour of $\al$.  This uses $v-2k$ colours to colour all $k$-subsets not contained in $\{1,\ldots,2k\}$.  The subsets not already coloured induce a copy of $\kne{2k}k$; since this graph is bipartite we can colour it with two colours.  Thus we have coloured $K_{v:k}$ with $v-2k+2$ colours.  Lov\'asz proved in \cite{lovkne} that this upper bound is the correct value.

There are at least three reasons why Lov\'asz's result is interesting.
First if $v=3k-1$, then $\knvk$ is triangle-free and has chromatic number $k+1$.
Hence we have an explicit construction of triangle-free graphs with large chromatic number.  By choosing $v$ and $k$ more carefully, we actually obtain graphs with large chromatic number and no short odd cycles.  (See \cite{cggr} for more details.)  Second, the fractional chromatic number of $\knvk$ is known to be $v/k$ and so the Kneser graphs provide examples of graphs whose fractional chromatic number is much lower than their chromatic number.  Third, Lov\'asz's proof that the chromatic number of $\knvk$ is $v-2k+2$ uses the Borsuk-Ulam theorem from topology in an essential way.  Although other proofs are known now, they all are based on results that are at least morally equivalent to the Borsuk-Ulam theorem.

Next we consider the connection between the ordinary Kneser graphs and the $q$-Kneser graphs.
To clarify this, we need the so-called $q$-binomial coefficients.  Choose a positive integer $q$ and, for an integer $n$, define
\[
[n] := \frac{q^n-1}{q-1}.
\]
We define the \textsl{$q$-factorial function} $[n]!$ inductively by $[0]!=1$ and
\[
[n+1]! =[n+1][n]!.
\]
We define the \textsl{$q$-binomial coefficient} $\gbin{v}k$ by
\[
\gbin{v}{k} =\frac{[v]!}{[k]![v-k]!}.
\]
This is also known as the Gaussian binomial coefficient, and is sometimes
written $\gbin{v}{k}_q$.
The $q$-binomial coefficients play the same role in the enumeration of subspaces that the usual binomial coefficients play in the enumeration of subsets.  If $q=1$, then we define $[n]$ to be $n$ and it then follows that $[n]! =n!$ and
\[
\gbin{v}{k} =\binom{v}k.
\]
If $q$ is a prime power, then $\gbin{v}k$ is equal to the number of subspaces of dimension $k$ in a vector space of dimension $v$ over the field of order $q$.  We recall, from \cite[p.~239]{blue}
for example, that in a $v$-dimensional vector space over $GF(q)$, the number of
$\ell$-dimensional subspaces that intersect a given $k$-dimensional subspace in a given subspace of dimension $j$ is
\[
q^{(\ell-j)(k-j)}\gbin{v-k}{\ell-j}.
\]

\begin{lemma}
The $q$-Kneser graph $\qkvk$ has $\gbin{v}k$ vertices and is regular with valency
$q^{k^2}\gbin{v-k}k$.
\end{lemma}

\proof
Only the valency is in question.  Suppose $|\fld|=q$ and $\al$ is a subspace of dimension
$k$ in $\fld^v$.  The subspaces of dimension $2k$ that contain $\al$ partition the set of
$k$-dimensional subspaces that meet $\al$ in the zero subspace.   This partition has
$\gbin{v-k}k$ components.  The number of $k$-subspaces in a space of dimension $2k$
that meet a given $k$-dimensional subspace in the zero subspace is $q^{k^2}$.\qed

When $q=1$, these expressions for the number of vertices and valency reduce to the corresponding values for the Kneser graph $\knvk$, respectively $\binom{v}k$ and $\binom{v-k}k$.  Many other parameters of the $q$-Kneser graphs, for example the eigenvalues of the adjacency matrix and their multiplicities, are given by expressions which involve $q$-binomial coefficients and which reduce to the corresponding value for the ordinary Kneser graphs when we set $q$ equal to $1$. In particular the fractional chromatic number of the ordinary Kneser graphs is $v/k$ while the fractional chromatic number of the $q$-Kneser graphs is $[v]/[k]$. (See \cite{cggr} for more details on fractional
chromatic number.)

Because of the above connections, and because the colouring problem for the ordinary Kneser graphs is so interesting, it is reasonable to study colouring problems for the $q$-Kneser graphs.
There is evidence that the relation between these two problems is complex.
Lov\'asz's result makes use of the fact that if $\seq\al1r$ are vertices in $\knvk$, then their set of common neighbours consists of the $k$-sets in the complement of the union
\[
\bigcup_{i=1}^r \al_i,
\]
which is a Kneser graph on a smaller set.  If $\seq\al1r$ is a set of vertices in $\qkvk$, then their set of common neighbours does not depend only on the join of these vertices.  Hence it is not easy to see how topological methods can be applied to colouring $q$-Kneser graphs.  It could be argued
that this adds interest to the $q$-colouring problem---it is not unreasonable to hope that
real progress on colouring $q$-Kneser graphs will yield insights concerning the case $q=1$.
However, in this paper we show that $\chi(\qkne{v}2)=[v-1]$ for $v>4$, and thus $\qkne52$ has chromatic number $[4]$.  But the ordinary Kneser graph $K_{5:2}$ is the Petersen graph, whose chromatic number is three.  So putting $q=1$ in this formula leads to the wrong answer.

\section{Independent Sets}

We will require information about the independent sets of maximum size in the ordinary Kneser graphs and the $q$-Kneser graphs, and so we summarise this here.

If $v \geq 2k$, the $k$-sets of $\{1,\ldots,v\}$ that contain a given point $i$ form an independent set in $\knvk$ with size $\binom{v-1}{k-1}$ and this is the maximum possible size.  If $v > 2k$ these are the only independent sets of this size but if $v=2k$, there are many others---partition the $k$-sets into $\binom{v-1}{k-1}$ pairs, and choose one $k$-set from each pair.

The Erd\H{o}s-Ko-Rado theorem asserts that if $v\ge(k-t+1)(t+1)$
and $\cF$ is a collection of $k$-subsets of a $v$-set such that
any two $k$-subsets in $\cF$ have at least $t$ elements in common,
then
\[
|\cF| \le\binom{v-t}{k-t}
\]
and, if $v>(k-t+1)(t+1)$ and equality holds, then $\cF$ consists of the $k$-subsets that contain
a given $t$-subset.  (This result was originally proved by Erd\H{o}s, Ko and Rado under the
assumption that $v$ was large enough relative to $k$, and it was first proved in the form stated
by Wilson \cite{rmwekr}.)

The optimal colourings of $\knvk$ we described earlier consist of $v-2k$ independent sets, each contained in an independent set of maximum size, together with a bipartition of $K_{2k:k}$.
To prove that $\chi(\knvk)=v-2k+2$ when $v>2k$, it would suffice to show the following: in any colouring of $\knvk$ with $v-2k+2$ colours, there is at least one colour class consisting of
$k$-subsets with a common point.  (The topological arguments offer no help here; they show
that there is no colouring using fewer than $v-2k+2$ colours.)

We turn to the $q$-Kneser graphs.  Frankl and Wilson \cite{fwqkne} proved that if $\cF$ is a collection of $k$-subspaces of $\fld^v$ such that any two elements of $\cF$ intersect in a subspace of  dimension at least $t$, then
\[
|\cF| \le\max\left\{\gbin{v-t}{k-t},\gbin{2k-t}k\right\}.
\]
When $v\ge2k$ and $t=1$, this implies that an independent set in $\qkne{v}k$ has size at
most $\gbin{v-1}{k-1}$.

What of the sets that meet this bound?  Here we must pay careful attention to the wording in \cite{fwqkne}.  It is asserted that when $v>2k$, if $\cF$ has maximum size then it consists of the $k$-spaces that contain a specified $t$-space.  No proof is offered, instead there is a claim that this result follows easily from the work in \cite{fwqkne} and the results in a second paper.  The difficulty is that the second paper is cited incorrectly; it is true that the characterisation follows readily from the work in \cite{fwqkne} and Wilson's earlier paper \cite{rmwekr}.  For the case $t=1$, a short proof is also offered in \cite{cgmn}.

When $v=2k$, the set of $k$-spaces contained in a given subspace of dimension $2k-1$
form an intersecting family with the maximum possible size.  Frankl and Wilson state that
that they could not prove there are only two types of optimal families when $t\ge2$.
Thus they assert that they can prove this when $t=1$, but it is not clear what proof they had in mind.  (M. W. Newman and the second author also have a proof of this now.)

\medbreak In our later work, we require information about minimal
sets of points in projective space that meet all subspaces of a
given dimension.  The basic result is the following.

\begin{theorem}
Let $S$ be a set of points in the  projective space $PG(v-1,q)$ such that every subspace of
projective dimension $k-1$ contains a point of $S$.  Then $|S|\ge[v-k+1]$, and if
equality holds then $S$ consists of the points from a subspace.\qed
\end{theorem}

The result was first proved by Bose and Burton \cite{bobu}.  As we
will need the case where we have a set of points $S$ that meet all
lines, we offer a proof for this case.  If $x$ is a point not in
$S$ then each line through $x$ must contain a point from $S$, and
so $|S|$ is bounded below by the number of lines on $x$.  If
equality holds then a line through $x$ meets $S$ in at most one
point.  Therefore if $\ell$ is a line that contains two points of
$S$, then all points on $\ell$ must lie in $S$.  Consequently $S$
is a subspace.

\section{Homomorphisms}
\label{sec:homkne}

If $X$ and $Y$ are two graphs, a homomorphism from $X$ to $Y$ is a map $f$ from $V(X)$ to $V(Y)$ such that if $u$ and $v$ are adjacent vertices in $X$, then $f(u)$ and $f(v)$ are adjacent in $Y$.  Since the graphs in this paper do not have loops, if $y\in V(Y)$, then the preimage $f\inv(y)$ is an independent set in $X$.  Further, $X$ can be coloured with $r$ colours if and only if there is a homomorphism $f:X\to K_r$.  It also follows that if there is a homomorphism
$f:X\to Y$, then $\chi(X)\le\chi(Y)$.  Accordingly homomorphisms provide a useful tool for
working on colouring problems.

The following homomorphisms between Kneser graphs are known.  First $K_{v:k}$ is an induced subgraph of $K_{v+1:k}$, we call this embedding the {\sl extension map}.  Next, if $t$ is a positive integer, then $K_{v:k}$ is an induced subgraph of $K_{tv:tk}$.  We
call this the {\sl multiplication map}.  Finally, Stahl \cite{stahl76} discovered a homomorphism from $K_{v:k}$ to $K_{v-2:k-1}$.  Given the existence of this map, we see that
\[
\chi(K_{v:k}) \le \chi(K_{v-2:k-1}),
\]
which implies that $\chi(K_{v:k})\le v-2k+2$.  Hence we call Stahl's map the {\sl colouring map}.  Stahl has conjectured that there is a homomorphism from $K_{v:k}$ to $K_{w,\ell}$ if and only if there is a homomorphism from $K_{v:k}$ to $K_{w,\ell}$  that is a composition of
extension, multiplication and colouring maps.

We turn to the $q$-Kneser graphs.  We have the following homomorphisms:
\begin{enumerate}[(a)]
\item
  The extension map, embedding $qK_{v:k}$ in $qK_{v+1:k}$.
\item
  Since the field of order $q$ is a subfield of the field of order $q^r$, we have a subfield map
  $qK_{v:k}\to q^rK_{v:k}$.
\item
  A $k$-space in a $v$-dimensional vector space over $GF(q^r)$ can be viewed as a subspace of
  dimension $rk$ in a space of dimension $rv$ over $GF(q)$.  This leads to a $q$-analog of the
  multiplication map, embedding $q^rK_{v:k}$ as an induced subgraph of $qK_{rv:rk}$.
\item
  Each $k$-subspace is the row space of a unique $k\times v$ matrix in reduced row echelon form.
  The subspace spanned by the last $k-1$ rows of this matrix is a $(k-1)$-subspace of a
  $(v-1)$-dimensional space.  Hence we have a homomorphism from $qK_{v:k}$ to
  $qK_{v-1:k-1}$.
\item
  Finally $qK_{v:k}$ is an induced subgraph of $K_{[v]:[k]}$.
\end{enumerate}
In Cases (c) and (e) above, the induced subgraph has the same fractional
chromatic number as the target graph.

However there is no homomorphism from $qK_{5:2}$ to $qK_{3:1}$, because
\[
\chi(qK_{3:1}) = q^2+q+1
\]
while
\[
\chi(qK_{5:2}) \ge \frac{\gbin{5}{2}}{\al(qK_{5:2})}
  =\frac{\gbin{5}{2}}{\gbin{4}{1}}
  = \frac{[5]}{[2]}  =\frac{q^4+q^3+q^2+q+1}{q+1} > q^3+q.
\]
Hence there is no $q$-analog of Stahl's colouring homomorphism
from $qK_{v:k}$ to $qK_{v-1:k-2}$.

\section{Chromatic Number}
\label{sec:chi}

Now we consider the chromatic number of the $q$-Kneser graphs.  There are two obvious
families of independent sets in $\qkne{v}{k}$, namely the set of all $k$-spaces containing a given
1-dimensional subspace, and the set of all $k$ spaces contained in a
given subspace of dimension $2k-1$.  There are $\gbin{v-1}{k}$ $k$-spaces containing
a given 1-dimensional subspace, while a subspace of dimension $2k-1$ contains
$\gbin{2k-1}{k}$ $k$-spaces.

\begin{lemma}
\label{lem:bnds}
  If $v\ge 2k$ then
\[
  \chi(qK_{v:k}) \le [v-k+1].
\]
  If $v=2k$ then
\[
  \chi(qK_{v:k})\le q^k+q^{k-1}
\]
\end{lemma}

\proof
A subspace $U$ of dimension $v-k+1$ has non-trivial intersection with each $k$-space, so we can colour each $k$-space $S$ with any of the 1-dimensional subspaces in $S \cap U$. As $U$ contains $[v-k+1]$ such subspaces, this yields a colouring with $[v-k+1]$ colours.

If $v=2k$ we can do even better than this by choosing a subspace
$U$ of dimension $k+1$, and a subspace $T$ of dimension $k$ in
$U$. Now consider the 1-dimensional subspaces in $U$ that do not
lie in $T$, together with the subspaces of dimension $2k-1$ that
contain $T$ but not $U$. This gives a total of
\[
[k+1]-[k] + [k]-[k-1] = [k+1]-[k-1] =q^k+q^{k-1}
\]
points and subspaces.

For any $k$-space $S$, if $S \cap U \subseteq T$, then $S$ lies in
a $(2k-1)$-space that contains $T$ but not $U$, and otherwise $S$
contains a 1-dimensional subspace of $U$ that does not lie in $T$.
Therefore we can use the points and subspaces as colours and
obtain a colouring of $\qkne{2k}{k}$ with $q^k + q^{k-1}$ colours.
\qed

The bound $[v-k+1]$ on $\qkne vk$ also follows from the fourth homomorphism described in
\sref{homkne}.  This leads to a second description of the colouring: if we represent each $k$-space by a $k\times v$ matrix in reduced row-echelon form, we can colour each subspace with
the first row.

It seems plausible to us that the upper bounds in \lref{bnds} provide the correct value of the
chromatic number in all cases.

\section{Covering Lines}
\label{sec:lines}

For the remainder of this paper we specialise to the situation $k=2$, where it proves
convenient to use the terminology of projective geometry. In this terminology, the 1-dimensional,
2-dimensional and 3-dimensional subspaces of $\fld^v$ are the points, lines and planes of
$PG(v-1,q)$. Two subspaces are called \textsl{incident} if one contains the other.
The $q$-Kneser graph $\qkne{v}{2}$ has the lines of $PG(v-1,q)$ as its vertices,
with two lines being adjacent if they are skew (have no point in common).

An independent set of size three
in $\qkne{v}2$ consists either of three concurrent lines, or three non-concurrent lines in the same plane.  It follows that an independent set of maximum size consists either of the lines on a point or the lines in a plane, and further that any independent set of $\qkne{v}{2}$ is contained in a maximum independent set of one of these types. Any colouring of $\qkne{v}{2}$ thus defines a collection of points and planes of $PG(v-1,q)$ such that every line of $PG(v-1,q)$ is incident with one of the points or one of the planes. We will call such a set of points and planes a \textsl{cover} of $PG(v-1,q)$ and say that a point or plane \textsl{covers} the lines with which it is incident.

If the colour classes of a colouring each contain more than $q+1$ vertices, then the colouring determines a unique cover. The converse is not quite true, in that  a cover does not determine a unique colouring of $\qkne{v}{2}$ because some lines may be incident with more than one element of a cover. However if a cover is minimal (under inclusion), all of the colourings it determines use the same number of colours, and so $\chi(\qkne{v}2)$ is equal to the minimum size of a cover.

\section{Projective 3-Space}
\label{sec:prtrs}

In this section we show that $\chi(\qkne42)=q^2+q$.

\begin{lemma}
\label{lem:colpts}
  Suppose $C$ is a cover of $PG(3,q)$ consisting of $r$ points and $s$ planes.  If $C$ contains
  $q+1$ collinear points, then $r+s\ge q^2+q+1$.  (Dually, if $C$ contains $q+1$ planes on one line,
  then $r+s\ge q^2+q+1$.)
\end{lemma}

\proof Let $\ell$ be a line all of whose points are in $C$. Then
these points cover
\[
(q+1)(q^2+q)+1 = q^3+2q^2+q+1
\]
lines and therefore there are $q^4$ remaining lines to be covered.

Each point of $C$ not in $\ell$ covers at most $q^2$ lines not
already covered by the points of $\ell$.  Each plane of $C$ meets
$\ell$ in at least one point, and so covers at most $q^2$ lines
not already covered by the points of $\ell$.  So we need in total
at least $q^2$ points and planes to cover the $q^4$ uncovered
lines.\qed

\begin{lemma}
\label{lem:rsgeq}
  Suppose $C$ is a cover of $PG(3,q)$ consisting of $r$ points and $s$
  planes.  If $r+s\le q^2+q$, then $r,s\ge q$.
\end{lemma}

\proof Let $x$ be a point not in $C$. There are $q^2+q+1$ lines on
$x$, and each point of $C$ covers at most one line on $x$. Since
$r+s\le q^2+q$, by assumption, we must have $r>0$. Similarly,
$s>0$.

Suppose, for a contradiction that, $r\le q-1$. Then at least
$q^2+2$ of the $q^2+q+1$ lines on $x$ are not covered by one of
the $r$ points. Hence they must be covered by one of the $s$
planes. The first plane on $x$ covers $q+1$ lines through $x$,
each additional plane on $x$ covers at most $q$ further lines.
Hence if there are $t$ of our $s$ planes on $x$, then
\[
(q+1)+(t-1)q \ge q^2+2
\]
and therefore
\[
t-1 \ge q-1+\frac{1}{q}.
\]
This implies that $t\ge q+1$.

Now count pairs $(x,H)$ where $x$ is a point not in our cover and
$H$ is a plane in the cover that contains $x$.  We find that
\[
s(q^2+q+1) \ge (q+1)[(q^3+q^2+q+1)-(q-1)]
\]
and hence
\[
s\ge \frac{(q^3+q^2+q-(q-2))(q+1)}{q^2+q+1}
  =q(q+1)-\frac{(q+1)(q-2)}{q^2+q+1}.
\]
Since $s$ is an integer, this implies that $s\ge q^2+q$.
Consequently $r+s \ge q^2+q+1$, which contradicts our initial
assumption. \qed

\begin{theorem}
  Suppose $C$ is a cover of $PG(3,q)$ with $r+s\leq{q^2+q}$ points and planes.  Then
  $C$ contains exactly $q^2+q$ points and planes, and, moreover, $q\mid r$ and dually $q \mid s$.
\end{theorem}

\proof
Suppose that $r=kq+x$ where $0\leq{x}<q$.  By \lref{rsgeq} we see that $s\ge q$, and therefore $r\leq{q^2}$ and so $k\leq{q}$.

Let $P$ be a plane that is not in the cover.  If $P$ does not contain a point of the cover, then
every line in $P$ would have to be covered by one of the planes, and therefore $s\ge q^2+q+1$ which is not possible.

Therefore $P$ contains points from the cover.  Suppose it contains at most $k$.  Since $k\leq{q}$, any $k$ points on $P$ cover at most $kq+1$ lines, so at least $(q^2+q+1)-(kq+1)=(q^2+q)-kq$ lines on $P$ are not covered by one of these $k$ points.  Each of these lines must be covered by one of the $s$ planes, whence $s\ge(q^2+q)-kq$ and $r+s\ge q^2+q+x$.  Therefore we may assume that any plane not in the cover contains at least $k+1$ of the points from the cover.

Each point of $PG(3,q)$ lies in exactly $q^2+q+1$ lines; since $r\le q^2+q$ it follows that
any point not in the cover lies on a line that contains no point from the cover.
Let $\ell$ be a line that does not contain points from the cover.  Since $r=kq+x$ where $x<q$, at least two of the $q+1$ planes on $\ell$, $P_1$ and $P_2$, will each contain fewer than $k+1$ points from the cover.  By the preceding paragraph, $P_1$ and $P_2$ must be in the cover.  To summarise, if $\ell$ is a line containing no points from the cover, then at least two planes on $\ell$ are in the cover.

The plane $P_1$ on $\ell$ contains at most $k$ points from the cover, so at least $(q^2+q)-kq$ lines on $P_1$ do not contain a point from the cover.  By the preceding paragraph each of these lines lies on a second plane from the cover, and we need at least $(q^2+q)-kq$ additional planes from the cover on these lines.  Consequently, $s\geq{(q^2+q)-kq}$ and $r+s\ge q^2+q+x$.
As $C$ satisfies $r+s\leq{q^2+q}$, we must have $x=0$. Therefore $r+s=q^2+q$ and, moreover, $q\mid r$.\qed

Since $\chi(\qkne42)$ is equal to the minimum size of a cover of
$PG(3,q)$, we have $\chi(\qkne42)=q^2+q$.

\section{Minimal Covers}

We have shown that $\chi(\qkne42)=q^2+q$, and given examples of
colourings which meet this bound.  In this section we completely
characterise the minimal covers of $PG(3,q)$.  Consider a cover
constructed as follows: Choose a plane $H$, a point $x$ on $H$,
and $s$ lines in $H$ on $x$, where $1 \le s \le q$.  The cover
then consists of the $q(q+1-s)$ points of $H$ not on these lines
and the $sq$ planes distinct from $H$ that contain one of the $s$
lines.  We will call a cover of this type a \textsl{standard
cover}.

\begin{lemma}
\label{lem:kpts}
Let $C$ be a cover of $PG(3,q)$ with $r=kq$ points and $s=q(q+1-k)$ planes.  If $P$
is a plane not in $C$, then it contains at least $k$ points from $C$; if equality holds then the
$k$ points are collinear.
\end{lemma}

\proof Since $P$ is not a plane in $C$, $P$ must contain some
points from $C$; otherwise, we'd need $q^2+q+1$ planes to cover
the lines on $P$. The $s=(q+1-k)q$ planes in $C$ cover at most
$(q+1-k)q$ lines on $P$ so at least $kq+1$ lines on $P$ remain. If
$t<q+1$, then $t$ points in a plane cover at most $tq+1$ lines,
with equality if and only if they are collinear. Consequently $P$
contains at least $k$ points from $C$, and if it contains exactly
$k$, they are collinear.\qed

\begin{lemma}
\label{lem:qplus1}
Let $C$ be a cover of $PG(3,q)$ with $r=kq$ points and $s=q(q+1-k)$ planes.  Suppose $P$ is a plane that contains at least $q+1$ points from the cover.  Then the points of $C$ cover at most $kq^3$ lines not in $P$.
\end{lemma}

\proof A point of $C$ in $P$ covers $q^2$ lines not in $P$.  If
$x$ is a point in $C$ not on $P$, at least $q+1$ lines incident
with $x$ are incident with a point in $C$ that lies on $P$.
Consequently $x$ covers at most $q^2$ lines that are not incident
with a point of $C$ on $P$.  Therefore the $kq$ points of $C$
cover at most $kq^3$ lines not on $P$.\qed

\begin{lemma}
\label{lem:aline}
If $C$ is a cover of $PG(3,q)$ with $r=kq$ points and $s=q(q+1-k)$ planes, then
each plane contains at least one line disjoint from $C$.
\end{lemma}

\proof
Let $P$ be a plane and assume by way of contradiction that each line on $P$ is incident
with a point from $C$.  Since all planes from $C$ meet $P$ in a line, each plane in $C$
contains at least one point from $C$.  Therefore on each plane of $C$ there are at least $q+1$
lines that are incident with a point from $C$, and consequently each plane in $C$ covers at
most $q^2$ of the lines that do not contain points from $C$.  As there are
\[
(q^2+1)(q^2+q+1) =q^4+q^3+2q^2+q+1
\]
lines in total, the number of lines covered by the points in $C$ is at least
\[
q^4+q^3+2q^2+q+1 -sq^2 =q^4+q^3+2q^2+q+1-(q+1-k)q^3 =kq^3+2q^2+q+1.
\]
We know that $q^2+q+1$ of these lines lie in $P$, the remaining lines, of which there
are at least $kq^3+q^2$, must intersect $P$ in a point.

Since every line in $P$ is incident with a point from $C$, there are at least $q+1$ points from $C$
on $P$.  By \lref{qplus1}, the $kq$ points of $C$ cover at most $kq^3$ lines not on $P$, a contradiction.\qed

\begin{lemma}
\label{lem:sqstd}
Let $C$ be a cover of $PG(3,q)$ with $q^2$ points and $q$ planes.  Then $C$ is standard.
\end{lemma}

\proof
Assume $C$ contains $q^2$ points and $q$ planes.  We will show first that there is a plane containing at least $q+1$ points from $C$.

Let $\ell_1$ be a line not incident with a point in $C$.  By \lref{colpts} there is a plane $H$
 on $\ell_1$ that is not in $C$ and by \lref{kpts}, there are at least $q$ points from $C$ on $H$.  If there are exactly $q$ points then they lie on a line $\ell_2$, and any plane containing $\ell_2$ and a point in $C$ not on $\ell_2$ contains $q+1$ points from $C$.  Otherwise $H$ contains at least $q+1$ points from $C$.

We next show that no plane of $C$ contains a point of $C$, and that the $q$ planes
of $C$ lie on a common line.

Let $P$ be a plane that contains at least $q+1$ points from $C$.  By \lref{qplus1}, the $q^2$
points in $C$ cover at most $q^4$ lines not in $P$.  By \lref{aline}, there is a line on $P$ that contains no point of $C$ and so at most $q^2+q$ lines on $P$ are incident with points of $C$.  Hence the number of lines incident with the $q^2$ points in $C$ is at most $q^4+q^2+q$.

Since any two planes have a line in common, the $q$ planes in $C$ cover at most $q^3+q^2+1$ lines.  The total number of lines is
\[
q^4+q^3+2q^2+q+1 =(q^4+q^2+q)+(q^3+q^2+1),
\]
whence the $q^2$ points in $C$ must cover exactly $q^4+q^2+q$ lines and the $q$ planes must cover exactly $q^3+q^2+1$.  We also see that the set of lines covered by the points of $C$ is
disjoint from the set of lines covered by the planes, and consequently no point of $C$ can lie
in a plane of $C$.

Further, since the $q$ planes cover exactly $q^3+q^2+1$ lines, the $q$ planes must lie on a line
$\ell$.

Let $Q$ be the unique plane on $\ell$ not in the cover.  Then the $q^2$ points of our cover must
lie on $Q$, and hence the points of the cover are the points of $Q\diff\ell$.\qed

\begin{lemma}
\label{lem:2planes}
Let $C$ be a cover of $PG(3,q)$ with $r=kq$ points and $s=q(q+1-k)$ planes and let
$P$ be a plane not in $C$ that contains at least $k+1$ points from $C$.  Then any line on $P$ not incident with a point from $C$ lies on at least two planes from $C$.
\end{lemma}

\proof Let $\ell$ be a line on $P$ that is not incident with a
point from the cover.  Let
\[
H_1,\ldots,H_{q+1}
\]
denote the $q+1$ planes on $\ell$, where $H_1=P$.  Since $\ell$
contains no point of the cover, these planes partition the $kq$
points of the cover.  By \lref{kpts}, each plane not in $C$
contains at least $k$ points from $C$.  Since $P$ contains $k+1$
points from the cover, it follows that at least two of the planes
on $\ell$ must lie in the cover.\qed

\begin{lemma}
\label{lem:LTP2}
Let $C$ be a cover of $PG(3,q)$ with $r=kq$ points and $s=q(q+1-k)$ planes.
Suppose there is a plane $P$ that contains at least $q+1$ points from the cover and a point $y$
that lies on at least $q+1$ planes.  Then no plane in the cover contains a point from the cover.
\end{lemma}

\proof
Let $P$ be a plane that contains at least $q+1$ points from $C$.  By \lref{qplus1}, our $r=kq$ points cover at most $kq^3$ lines not in $P$.  Since there is a line in $P$ that contains no points from $C$, our $r$ points cover at most $kq^3+q^2+q$ lines.

Dually, the number of lines covered by the $s$ planes in $C$ is at most
\[
sq^2+q^2+q= (q+1-k)q^3+q^2+q.
\]

Suppose, for a contradiction, that some plane in $C$ contains a point from $C$.  Then there are
$q+1$ lines that are covered by both by a point in $C$ and a plane from $C$.  Hence the
number of lines covered by the points and planes of $C$ is at most
\[
(kq^3+q^2+q)+((q+1-k)q^3+q^2+q)-(q+1) =q^4+q^3+2q^2+q-1
\]
Since there are $q^4+q^3+2q^2+q+1$ lines altogether, this provides our contradiction.\qed

\begin{theorem}
A cover of $\qkne42$ with $q^2+q$ points and planes is standard.
\end{theorem}

\proof
Let $C$ be a cover of $PG(3,q)$ with $q^2+q$ points and planes.  We may assume
that there are $r=kq$ points and $s=(q+1-k)q$ planes.  By  \lref{sqstd} and
duality, we may assume that $2\le k\le q-1$.

As a first step, we show that there is a line that contains no point from $C$ and lies on exactly one
plane from $C$.  Let $m$ be a line that contains no point from $C$.  There are $q+1$ planes on $m$ and $kq$ points in the cover, so there is a plane $H$ on $m$ that contains fewer than $k$ points.  By \lref{kpts} we see that $H$ lies in the cover.  At most $(k-1)q+1$ lines on $H$ are
incident with points of $C$ and therefore there are at least $q(q+2-k)$ lines in $H$ not incident with a point from $C$.  Since there are only $q(q+1-k)$ planes in $C$, there is a line $\ell$ in $H$ which contains no point from $C$ and which is not contained in a second plane from $C$.

Next we show that there is a plane that contains at least $q+2$ points from $C$.

Let $\seq H1q$ denote the planes on $\ell$ other than $H$.  These $q$ planes do not belong to $C$ and therefore by \lref{kpts}, there are at least $k$ points from $C$ on each of them.  Since these planes partition the points of $C$ into $q$ classes, each plane contains exactly $k$ points from $C$ and, by \lref{kpts}, each set of $k$ points lies on a line.  Denote the line on $H_i$ by $m_i$.
The $k$ points on $H_i$ cover exactly $kq+1$ lines on $H_i$; the remaining $(q+1-k)q$ lines on $H_i$ are covered by planes of $C$.  Since there are exactly $(q+1-k)q$ planes in $C$,
each line of $H_i$ that is not covered by a point of $C$ is contained in exactly one plane from $C$.  Note that $H$ contains no points of $C$.

The plane $H_1$ contains $\ell$ and therefore $m_1$ intersects $\ell$ in a point $x$. The lines
other than $m_1$ on $x$ in $H_1$ are covered by planes of $C$, and so there are $q$ planes
from $C$ on $x$.  For $i=2,\ldots,q$ these planes intersect $H_i$ in $q$ distinct lines
through $x$, and these lines do not contain points of $C$.  Therefore each of the lines $\seq m2q$ intersects $\ell$ in $x$.

Let $P$ be the plane determined by $m_1$ and $m_2$.  The lines on $P$ incident with $x$
are $m_1$, $m_2$, the intersection of $P$ with $H$ and the intersection of $P$ with
$\seq H3q$.  The planes in $C$ intersect $H_1$ in lines that contain no points of $C$,
but $P\cap H_1=m_1$ which does contain points from $C$.  Therefore $P$ is not in $C$.

Since $P$ contains $2k$ points from $C$, by \lref{2planes} any line on $P$ not incident
with a point from $C$ lies on at least two planes from $C$.  As there are $q$ planes from $C$ on $x$, at most $q/2$ lines on $P$ incident with $x$ do not contain points from $C$.  Consequently
at least $(q+2)/2$ lines on $P$ incident with $x$ contain points from $C$.  Referring to our
listing above of the lines on $x$ in $P$, we see that $m_1$ and $m_2$ contain $k$ points from
$C$.  As $H$ contains no points of $C$, the line $P\cap H$ is disjoint from $C$.  If $P\cap H_i$
contains a point from $C$ then $P\cap H_i$ is $m_i$, because this is the only line on $x$ in $H_i$
that contains points from $C$.  Then $P\cap H_i$ contains $k$ points from $C$.
Since $k\ge2$, it follows that the number of points from $C$ on $P$ is at least
\[
k\frac{q+2}2 \ge q+2.
\]
So we have shown that there is a plane that contains at least $q+2$ points from $C$; the dual of
our argument shows that there is a point $y$ on at least $q+2$ planes from $C$.

By \lref{LTP2}, no plane in $C$ contains a point from $C$.  The $(q+1-k)q$ planes in $C$ each meet $P$ in a line, and so by \lref{colpts} there are at least $q+1-k$ lines in $P$ that contain no point of $C$.  Thus there are at most $q^2+k$ lines in $P$ that do contain points of $C$.
Each point of $C$ in $P$ covers $q^2$ lines not in $P$.  Since there are at least $q+2$ points of $C$ in $P$, each point of $C$ not in $P$ covers at most
\[
(q^2+q+1)-(q+2) =q^2-1
\]
lines not covered by points of $C$ in $P$.  So the number of lines covered by the points of $C$ is
at most
\[
kq^3+q^2+k,
\]
and if equality holds, all points of $C$ lie in $P$ and there are exactly $q+1-k$ lines in $P$
disjoint from $C$, each of which lies in $q$ planes from $C$.

Dually, the number of lines covered by the $(q+1-k)q$ planes of $C$ is at most
\[
(q+1-k)q^3+q^2+(q+1-k)
\]
and, if equality holds, these planes have a common point $y$ and there are exactly
$k$ lines incident with $y$ that do not lie on a plane from $C$.  Since the total number
of lines is
\[
q^4+q^3+2q^2+q+1 =(kq^3+q^2+k)+((q+1-k)q^3+q^2+q+1-k),
\]
our last two inequalities must be tight.  Therefore all points of $C$ lie in $P$ and all the planes of
$C$ contain $y$.  Hence $C$ is a standard cover.\qed

\section{Higher Dimensions}

In this section we determine the chromatic number of $\chi(\qkne{v}{2})$ for $v \geq 5$ by
determining the minimum number of points and planes needed to cover the lines
of $PG(v-1,q)$.

\begin{theorem}
Let $C$ be a cover of $PG(v-1,q)$, where $v \ge 5$ with $r$ points and $s$ planes
such that $r+s\le[v-1]$.  Then $r\ge[v-1]$.
\end{theorem}

\proof Suppose, for a contradiction, that $r<[v-1]$, and define
$\de:= [v-1]-r$. We determine a lower bound on the number of lines
that do not contain a point of $C$, by counting the flags
$(x,\ell_x)$ where $x$ is a point not in $C$ and $\ell_x$ is a
line on $x$ containing no points of $C$. Each point $x$ not in $C$
lies on $[v-1]$ lines, and therefore there are at least
$[v-1]-r=\de$ lines through $x$ that contain no points of $C$.  As
there are $[v]-r$ points not in $C$, the number of flags is at
least
\[
([v]-r)\de =([v]-[v-1]+[v-1]-r)\de =(q^{v-1}+\de)\de.
\]
Since each line containing no points of $C$ lies in exactly $q+1$ flags, it follows that the
number of such lines is at least
\[
\frac{(q^{v-1}+\de)\de}{q+1}.
\]

Each of these lines must be contained in one of
the $s$ planes, and a plane contains exactly $q^2+q+1$ lines.  Therefore
\begin{equation}
\label{sqd}
s \ge\frac{(q^{v-1}+\de)\de}{(q+1)(q^2+q+1)}.
\end{equation}

Since $r+s\le[v-1]$, we have $s\le\de$ so
\[
\frac{(q^{v-1}+\de)\de}{(q+1)(q^2+q+1)} \le\de,
\]
from which we have
\begin{equation}
\label{deqq}
\de \le(q+1)(q^2+q+1)-q^{v-1}.
\end{equation}

Observe that
\[
q^4-(q+1)(q^2+q+1) =q^4-q^3-2q^2-2q-1 =q(q(q^2-q-2)-2)-1
\]
and therefore
\[
q^{v-1}-(q+1)(q^2+q+1) =(q^{v-1}-q^4)+q(q(q^2-q-2)-2)-1.
\]
If $q>2$ then $q^2-q-2>0$ and so the right side is positive.  If
$q=2$ then the right side is equal to
\[
(2^{v-1}-16)-5,
\]
which is positive if $v>5$. Consequently, we conclude that the
right hand side of (2) is negative in these cases, which is a
contradiction. Therefore, $r=[v-1]$ if $q>2$ or if $v>5$ and
$q=2$.

\medbreak
Finally we consider the case where $v=5$ and $q=2$. Let
$x$ be a point not in $C$. Since $r<[4]=15$, by assumption, at
least one of the 15 lines on $x$ must be covered by one of the $s$
planes in the cover. Consequently, $x$ must lie on one of the $s$
planes in the cover. Since, $r+s\leq{[4]}=15$, we have
$r\leq{15-s}$ so at least
\begin{equation}
\label{31s}
31-(15-s)=16+s
\end{equation}
points don't lie in $C$, and must lie on one of the $s$ planes in
the cover. Since planes contain 7 points, we must have
$7s\geq{16+s}$ so $s\geq{3}$.

Suppose $s>3$.  Then $r\leq{11}$, so at least four of the 15 lines
on $x$ must lie on the $s$ planes. Since a plane on $x$ covers
three lines on $x$, we must have that $x$ lies on at least two
planes in the cover.  Consequently, ${7s}/{2}\geq{16+s}$, which
implies that $s\geq{7}$. Since $s\leq{\de}\leq{5}$ by (2), this is
a contradiction.

Therefore, $s=3$ so at least 19 points must lie on the three
planes in the cover by \eqref{31s}.  However in $2K_{5:2}$ distinct
planes intersect, so three planes cannot contain 19 distinct
points.  We have the desired contradiction, so $r=[4]=15$ when $v=5$
and $q=2$.\qed

\medbreak
\noindent{\sc Acknowledgement:}  The work in this paper has benefited from a number of
discussions with Ada Chan.

\end{document}